\documentclass[a4paper,12pt]{article}
\title{Refinement of asymptotic behavior of the eigenvalues for the linearized Liouville-Gel'fand Problem}
\date{}
\author{Hiroshi Ohtsuka\thanks{Faculty of Mathematics and Physics, Institute of Science and Engineering, Kanazawa University, Kakuma-machi, Kanazawa-shi, Ishikawa, 920-1192, Japan \\ e-mail \textsf{ohtsuka@se.kanazawa-u.ac.jp}}, \quad Tomohiko Sato\thanks{Department of Liberal Arts and Basic Sciences, College of Industrial Technology, Nihon University, 2-11-1 Shin-ei, Narashino, Chiba 275-8576, Japan \\ e-mail \textsf{sato.tomohiko@nihon-u.ac.jp}}}
\usepackage{amssymb}
\usepackage{amsmath}
\usepackage{color}
\usepackage{layout}
\usepackage{times}
\usepackage{enumerate}
\usepackage{mathtools}
\mathtoolsset{showonlyrefs=true}
\newtheorem{thm}{Theorem}[section]
\newtheorem{lem}{Lemma}[section]
\newtheorem{prop}{Proposition}[section]

\newtheorem{defi}{Definition}[section]

\begin{document}
\maketitle
\begin{abstract}
We determine the second term of the asymptotic expansions for the first $m$ eigenvalues and eigenfunctions of the linearized Liouville-Gel'fand problem associated to solutions which blow-up at $m$ points. Our problem is the case with an inhomogeneous coefficient in two dimensional domain and we extend the previous studies for the problem with a homogeneous coefficient. We also discuss  in detail the required regularity of the coefficient necessary to get the conclusion.
\mbox{} 
\end{abstract}
\section{Introduction and Preliminary}
Let $\Omega \subset \mathbb{R}^2$ be a bounded domain with smooth boundary $\partial \Omega$. We consider the Liouville-Gel'fand problem with Dirichlet boundary condition, 
\begin{equation}\label{f1.1}
-\Delta v = \lambda V(x) e^v \quad \text{in} \ \Omega, \qquad v = 0 \quad \text{on} \  \partial \Omega 
\end{equation}
where $\lambda >0$ is a parameter and $V=V(x)$ is a positive, suitably smooth function.
\par In this paper, we consider the blow-up solution to (\refeq{f1.1}) at interior $m$-point in $\Omega$, and we will give 
a refinement of the asymptotic behavior of $n$-th eigenvalue, $1 \leq n \leq m$, for the linearized eigenvalue problem 
\begin{align}
&- \Delta w^n = \mu^n \lambda V(x) e^{v} w^n \quad \text{in} \ \Omega, \qquad w^n =0 \quad \text{in} \ \partial \Omega, \\
&\|w^n \|_{L^\infty (\Omega)} = \max_{\overline{\Omega}} w^n =1
\end{align}
of (\refeq{f1.1}).
We also study the characterization of the corresponding eigenfunction $w = w^n$ at $m$ blow-up points. Similar problem is already studied in \cite{ggo16}, but they treat only the case $V(x)\equiv 1$. We extend it to the case to $V(x)\in C^1(\overline{\Omega})$ in this paper. In order to do this, we require new asymptotic estimate (Lemma \ref{prop3.3}, described below) which is, previously, only obtained under the stronger condition.
\par Now we review background theorems before we state our results. After the pioneering work of Nagasaki-Suzuki\cite{ns90} for $V\equiv 1$, it is well-known that the Hamiltonian controls the formation of singular limit of the solution to the problem (\refeq{f1.1}). The following theorem, which is an extension of the results of  \cite{ns90} for $V\equiv 1$ to some general $V(x)$, studies the formulation of singular limit of the blow-up solution to the equation of mean field type with Dirichlet boundary condition, and the location of blow-up points controlled by the Hamiltonian. 
\begin{thm}[\cite{mw01}]\label{th-mw01}
Assume $0<V=V(x) \in C^{1}(\overline{\Omega})$. Let $v=v_k$, $k=1,2,\ldots ,$ be a sequence of solutions to 
\begin{equation}\label{f0.7}
-\Delta v = \sigma \frac{V(x) e^v}{\int_{\Omega}V(x) e^v} \quad \text{in} \ \Omega, \qquad v = 0 \quad \text{on} \  \partial \Omega. 
\end{equation}
with $\sigma = \sigma_k \to \sigma_0 \in (0, +\infty)$, and we assume $\displaystyle \lim_{k \to \infty}\|v_k\|_{L^\infty(\Omega)} = + \infty$. Then, passing to subsequence, it holds that $\sigma_0 = 8 \pi m$ for some positive integer $m$ and there exists a blow-up set defined by
\begin{equation}\label{f0.2}
\mathcal{S} =\{x_0 \in \overline{\Omega} \ | \ \exists \{x_k\} \subset \Omega \text{ such that } x_k \to x_0 \text{ and } v_k(x_k) \to + \infty \}
\end{equation}
which is composed of distinct $m$-point $\{x^\ast_1, \ldots , x^\ast_m\} \subset \Omega$. Moreover, we have the singular limit 
\begin{equation}\label{f0.3}
v_k \to 8 \pi \sum^{m}_{j=1} G(\cdot , x^{\ast}_{j}) 
\end{equation}
locally uniformly in $\overline{\Omega}\setminus \mathcal{S}$, and it holds that $(x^\ast_1, \ldots , x^\ast_m)$ is a critical point of the function 
\begin{equation}\label{f0.8}
H^m_V (x_1, \cdots, x_m) = \frac{1}{2}\sum^{m}_{i=1}R(x_i) + \sum_{1 \leq i < j \leq m}G(x_i, x_j) + \frac{1}{8\pi}\sum^{m}_{i=1}\log V(x_i)
\end{equation}
which is the Hamiltonian of the system of $m$-point vortices defined for \\ $(x_1 , \ldots , x_m) \in \Omega^m \setminus D$, where $G=G(x,y)$ is the Green function of 
\begin{equation}\label{f0.5}
- \Delta G(\cdot , y) = \delta_y \quad \text{in}\ \Omega, \qquad G(\cdot , y)=0 \quad  \text{on}\ \partial \Omega,
\end{equation}
$D=\{(x_1, \ldots , x_m) \in \Omega^m = \Omega \times \cdots \times \Omega \ |\ \exists (i,j)  \text{ such that } i \neq j \text{ and } x_i = x_j\}$ and $R=R(x)$ is the Robin function $R(x) = K(x,x)$, $K(x,y)= G(x,y) + \frac{1}{2 \pi} \log |x-y|$ where $K \in C^2\left((\overline{\Omega}\times \Omega) \cup (\Omega \times \overline{\Omega})\right)$.
\end{thm}
We note that (\refeq{f0.3}) holds in $C^{2, \alpha}_{loc}(\overline{\Omega} \setminus \mathcal{S})$ from the elliptic regularity. We also note that the regularity of $V$ is not stated clearly in \cite{mw01}, but it is now well known that $V\in C^1(\overline{\Omega})$ is sufficient to get the above result.
\par Assume $0<V=V(x) \in C^{1}(\overline{\Omega})$. By putting $\lambda = \sigma / \int_{\Omega}V(x) e^v$ into (\refeq{f0.7}), we have the considered problem (\refeq{f1.1}). Let $v=v_k$ be a sequence of solutions to (\refeq{f1.1}) with $\lambda = \lambda_k \downarrow 0$. Assume $\| v_k \|_{L^\infty(\Omega)}  \to + \infty$ as $k \to \infty$, and we consider 
\begin{align}
\Sigma_k = \int_{\Omega} \lambda_k V e^{v_k} \to 8 \pi m \label{f0.11}
\end{align}
for some positive integer $m$.
\par Now we consider the linearized eigenvalue problem
\begin{align}\label{f1.2}
&- \Delta w^n_k = \mu^n_k \lambda_k V(x) e^{v_k} w^n_k \quad \text{in} \ \Omega, \qquad w^n_k =0 \quad \text{in} \ \partial \Omega, \\
&\|w^n_k\|_{L^\infty (\Omega)} = \max_{\overline{\Omega}} w^n_k =1
\end{align}
where $\mu^n_k$, $n = 1, 2, \ldots$ is the $n$-th eigenvalue and $w^n_k=w^n_k(x)$ is the eigenfunction corresponding to $\mu^n_k$. The eigenvalues $\mu^n_k$ satisfy the following relation
\begin{align}
0 < \mu^1_k < \mu^2_k \leq \mu^3_k \leq \cdots \to + \infty, \label{f0.13}
\end{align}
and we assume the orthogonality in $H^1_0(\Omega)$:
\begin{align}\label{f1.3}
\int_{\Omega} \nabla w^i_k \cdot \nabla w^j_k =0, \qquad i \neq j.
\end{align}
From Theorem \ref{th-mw01}, it follows that $\mathcal{S} =\{x^\ast_1 , \ldots , x^\ast_m\} \subset \Omega$ and $\|v_k\|_{L^{\infty}(\omega)} = O(1)$ for every $\omega \subset \subset \overline{\Omega}\setminus \mathcal{S}$.
\par Take $0 < R \ll 1$ such that $B_{2R}(x^\ast_j) \subset \subset \Omega$, $1 \leq j \leq m$ and 
$B_R(x^\ast_i) \cap B_R(x^\ast_j) = \emptyset$, $i \neq j$. For $x^\ast_j \in \mathcal{S}$, $1 \leq j \leq m$, there exists $\{x_{j,k}\}_k \subset B_R(x^\ast_j)$ such that
\begin{align}\label{f0.12}
v_k(x_{j,k}) = \sup_{B_R(x_{j,k})} v_k \to + \infty, \qquad x_{j,k} \to x^\ast_j
\end{align}
as $k \to \infty$ by (\refeq{f0.2}).
\par We introduce a scaling parameter $\delta_{j,k}>0$ defined by
\begin{equation}\label{f2.1}
\lambda_k V(x_{j,k}) e^{v_k(x_{j,k})} \delta_{j,k}^2 =1.
\end{equation}
We set $w_k := w^n_k$, $\mu_k := \mu^n_k$, and let 
\begin{align}
\tilde{v}_{j,k}(\tilde{x})&=v_k(\delta_{j,k}\tilde{x} + x_{j,k}) - v_k(x_{j,k}), \\
\tilde{w}_{j,k}(\tilde{x})&=w_k(\delta_{j,k}\tilde{x} + x_{j,k}), \label{f2.2}\\
\displaystyle \tilde{V}_{j,k}(\tilde{x})&=\frac{V(\delta_{j,k}\tilde{x} + x_{j,k})}{V(x_{j,k})}
\end{align}
for $\tilde{x} \in B_{R/\delta_{j,k}}(0)$. We have $\delta_{j,k} \to 0$ by \cite{bm91}, and we also have $\tilde{V}_{j,k}(\tilde{x}) \to 1$ locally uniformly in $\tilde{x} \in \mathbb{R}^2$. Moreover, by using the local uniform estimate of Y.Y.Li type \cite{yyli99}, 
there exists a constant $C>0$ such that
\begin{align}
\left| v_k(x) - \log \frac{e^{v_k(x_{j,k})}}{\left(1 + \frac{1}{8}\lambda_k V(x_{j,k}) e^{v_k(x_{j,k})} |x- x_{j,k}|^2 \right)^2} \right| \leq C \label{f6.6}
\end{align}
for any $x \in B_R(x_{j,k})$, $1 \leq j \leq m$ and $k \gg 1$. By (\refeq{f2.1}) and (\refeq{f6.6}), there exist constants $d_j >0$, $j = 1, \ldots , m$ such that 
\begin{align}
\delta_{j,k} = d_j \lambda_k^{1/2} + o(\lambda_k^{1/2}). \label{f6.4}
\end{align}
\par Concerning $\tilde{v}_{j,k}$ of (\refeq{f2.2}), by (\refeq{f1.1}), (\refeq{f0.11}) and (\refeq{f2.1}) it holds that 
\begin{align}\label{f2.3}
&\displaystyle - \Delta \tilde{v}_{j,k} = \tilde{V}_{j,k} e^{\tilde{v}_{j,k}}, \quad \displaystyle \tilde{v}_{j,k} \leq \tilde{v}_{j,k}(0) = 0 \quad \text{in} \ B_{R/\delta_{j,k}}(0),  \\
&\int_{B_{R/\delta_{j,k}}(0)} \tilde{V}_{j,k} e^{\tilde{v}_{j,k}} = \int_{B_R(x_{j,k})} \lambda_k V e^{v_k} \to 8\pi.
\end{align}
By \cite{cl91} and the elliptic regularity, we have that
\begin{equation}\label{f2.4}
\tilde{v}_{j,k}(\tilde{x}) \to U(\tilde{x}) = \log \frac{1}{\left(1 + \frac{1}{8} |\tilde{x}|^2 \right)^2}  
\end{equation}
locally uniformly in $\tilde{x} \in \mathbb{R}^2$ and $U= U(\tilde{x})$ satisfies
\begin{equation}\label{f2.5}
- \Delta U = e^{U}, \quad U \leq U(0) = 0 \quad \text{in} \ \mathbb{R}^2, \qquad \int_{\mathbb{R}^2}e^{U} < + \infty.
\end{equation}
\par From relations (\refeq{f2.1}) and (\refeq{f6.4}), it follows that
\begin{align}
\log \lambda_k + \log V(x_{j,k}) + v_k(x_{j,k}) + \log \left\{d_j^2 \lambda_k + o(\lambda_k) \right\} = 0.
\end{align}
Therefore, we obtain
\begin{align}
v_k(x_{j,k}) = -2 \log \lambda_k - 2 \log d_j - \log V(x_{j,k}) + o(1). \label{f6.5}
\end{align}
\par By using the estimate (\refeq{f6.6}), we have 
\begin{align}
\left| \tilde{v}_{j,k}(x) - U(x) \right| \leq C \label{f6.7}
\end{align}
for $x \in B_{R/ \delta_{j,k}} (0)$. 
\par On the other hand, concerning $\tilde{w}_{j,k}$ of (\refeq{f2.2}), it holds that
\begin{align}\label{f2.6}
&\displaystyle - \Delta \tilde{w}_{j,k} = \mu_k \tilde{V}_{j,k} e^{\tilde{v}_{j,k}} \tilde{w}_{j,k}  \quad \text{in} \ B_{R/\delta_{j,k}}(0), \\
&\displaystyle \|\tilde{w}_{j,k}\|_{L^{\infty}(B_{R/\delta_{j,k}}(0))} \leq 1
\end{align}
by (\refeq{f1.2}) and (\refeq{f2.1}). From Lemma 2.6 of \cite{ss18}, it follows that
\begin{align}
w_k = o(1) \quad \text{in} \ C^1 \left(\overline{\Omega}\setminus \cup^{m}_{j=1}B_R(x^\ast_j) \right).
\end{align}
Now we define the concentration of eigenfunction $w^n_k$ of (\refeq{f1.2}) in the same manner to \cite{ggo16}.
\begin{defi}
We say that an eigenfunction $w^n_k$ concentrates at $\kappa \in \Omega$ if there exist a constant $C>0$ and $\{\kappa_{k}\}_k \subset \Omega$ such that $\kappa_{k} \to \kappa$ and
\begin{align}\label{f2.7}
|w^n_k(\kappa_{k})| \geq C \quad \text{for} \quad  k \gg 1. 
\end{align}
\end{defi}
We recall a result in \cite{ss18} as an extension of that in \cite{ggos14} which indicates a rough estimate of eigenvalue $\mu^n_k$ and convergences of $\tilde{w}_{j,k}$ and $w^n_k / \mu^n_k$, $1 \leq n \leq m$.
\begin{thm}[\cite{ss18}]\label{th1.1}
Assume $0<V=V(x) \in C^{1}(\overline{\Omega})$. Let $v_k$, $k=1,2, \ldots ,$ be a sequence of solutions to (\ref{f1.1}) for $\lambda=\lambda_k$ and let $\lambda_k \downarrow 0$. Assume $\|v_k\|_{L^\infty (\Omega)} \to + \infty$, and consider $\Sigma_k = \int_{\Omega} \lambda_k V e^{v_k} \to 8 \pi m$ and $\mathcal{S} = \{x^\ast_1, \ldots , x^\ast_m\}$ in Theorem \ref{th-mw01}. Then, it holds that $\mu^n_k \to 0$ for each $n=1, \ldots , m$, and there exists
\begin{align}
\boldsymbol{c}^n = (c^n_1, \ldots , c^n_m) \in \mathbb{R}^m \setminus \{\boldsymbol{0}\} \label{f3.1}
\end{align}
such that for each $j= 1, \ldots , m$, taking a subsequence, 
\begin{align}
&|c^n_j| \leq 1, \\
&\tilde{w}^n_{j,k} \to c^n_j \qquad \text{locally uniformly in} \ \mathbb{R}^2, \label{f3.2} \\
&\boldsymbol{c}^i \cdot \boldsymbol{c}^{i'} =0, \qquad i \neq i' \label{f3.3}
\end{align}
and
\begin{align}
\frac{w^n_k}{\mu^n_k} \to 8 \pi \sum^m_{j=1} c^n_j G(\cdot, x^\ast_j) \quad \text{uniformly in} \  \overline{\Omega}\setminus \cup^{m}_{j=1} B_R(x_{j,k}) \label{f3.4}
\end{align}
where $G=G(x,y)$ is the Green function defined in (\refeq{f0.5}).
\end{thm}
We note that (\refeq{f3.4}) holds for $C^2(\overline{\Omega}\setminus \cup^{m}_{j=1} B_R(x_{j,k}))$ from the elliptic regularity similar to the case of (\refeq{f0.3}) for $V \in C^1(\overline{\Omega})$.
From Lemma 2.8 of \cite{ss18} and Theorem \ref{th1.1}, it follows that 
\begin{align}\label{f3.6}
w^n_k \ \text{concentrates at} \  x^\ast_j \quad \text{and} \quad c^n_j \neq 0 
\end{align}
are equivalent.
\par The following theorem asserts the asymptotic behavior of eigenvalues $\mu^n_k$ for the problem (\refeq{f1.2}). This theorem is also an extension of \cite{ggos14}.
\begin{thm}[\cite{ss18}]\label{th-ss18-1} We assume $V(x) \in C^{2}(\overline{\Omega})$ in (\refeq{f1.1}), in addition to the assumption of Theorem \ref{th1.1}. Then it holds that
\begin{align}
&\mu^n_k = - \frac{1}{2 \log \lambda_k} + o\left(\frac{1}{\log \lambda_k} \right), && 1 \leq n \leq m, \label{f0.10}\\
&\mu^n_k = 1-48 \pi \eta^{2m-(n-m)+1}\lambda_k + o(\lambda_k), && m+1 \leq n \leq 3m, \label{f0.10-2} \\
&\mu^n_k >1, && n \geq 3m+1, 
\end{align}
for $k \gg 1$. Here, $\eta^n$, $1 \leq n \leq 2m$ is the $n$-th eigenvalue of \\ $D[(\mathrm{Hess} H^m_V)(x^\ast_1, \ldots , x^\ast_m)]D$, where $H^m_V$ is the Hamiltonian defined in (\refeq{f0.8}) and $D$ is a diagonal matrix $D = \mathrm{diag}[d_1, d_1, d_2, d_2, \ldots , d_m, d_m]$ with $d_j >0$ appearing in (\refeq{f6.4}).
\end{thm}
%
%
%
%
We have the behavior (\refeq{f0.10}) under the same assumption of Theorem \ref{th1.1}, that is, $V \in C^1(\overline{\Omega})$, and we require the regularity $V(x) \in C^{2}(\overline{\Omega})$ to prove the behavior (\refeq{f0.10-2}). As a main result in this paper, we make refinement of the asymptotic behavior (\refeq{f0.10}) in Theorem \ref{th-ss18-1}. This theorem is natural extension of \cite{ggo16} considering the case of $V \equiv 1$ in (\refeq{f1.1}), and will be proven under $V(x) \in C^{1}(\overline{\Omega})$.
\begin{thm}\label{th1.6}
Under the same assumption of Theorem \ref{th1.1}, we have the following (i) and (ii):
\begin{enumerate}[(i)]
\item For each $n=1. \ldots , m$, the vector $\boldsymbol{c}^n = (c^n_1, \ldots , c^n_m)$ in Theorem \ref{th1.1} is the eigenvector corresponding to the $n$-th eigenvalue of the $m \times m$ matrix $(h_{ij})$ defined by
\begin{align}\label{f4.1}
h_{ij} =
\begin{cases}
\displaystyle R(x^\ast_i) + 2\sum^{}_{1 \leq h \leq m, \ h \neq i}G(x^\ast_h, x^\ast_i) + \frac{1}{4\pi}\log V(x^\ast_i), \qquad i=j,\\
-G(x^\ast_i, x^\ast_j), \qquad i \neq j.
\end{cases}
\end{align}
\item It holds that 
\begin{align}
&\mu^n_k \\
&= -\frac{1}{2 \log \lambda_k} + \left(2 \pi \Lambda^n - \frac{3\log 2 -1}{2} \right) \cdot \frac{1}{(\log \lambda_k)^2} + o\left( \frac{1}{(\log \lambda_k)^2} \right) \label{f5.2}
\end{align}
for $n = 1 ,\ldots , m$, $k \gg 1$, where $\Lambda^n$ is the $n$-th eigenvalue of the matrix $(h_{ij})$ defined in (\refeq{f4.1}), and we consider $\Lambda^1 \leq \cdots \leq \Lambda^m$.
\end{enumerate}
\end{thm}
%
%
\par The following theorem is another result regarding a characterization of $\boldsymbol{c}^n$ in (\refeq{f3.1}) concerning the matrix $(h_{ij})$.
\begin{thm}\label{th1.2}
Under the same assumption of Theorem \ref{th1.1}, for each $n=1. \ldots , m$ and for some subsequence of $\{w^n_k\}_k$, it holds that
\begin{align}
\tilde{w}^n_{j,k}(x) = w^n_k(x_{j,k}) +\mu^n_k c^n_j U(x) +o(\mu^n_k), \quad 1 \leq j \leq m \label{f4.2}
\end{align}
locally uniformly in $x \in \mathbb{R}^2$, where $U(x)$ is defined in (\refeq{f2.4}).
\end{thm}
We note that (\refeq{f4.2}) holds in $C^{2,\alpha}_{loc}(\mathbb{R}^2)$ by the elliptic regularity.
\par We have the following proposition by Theorem \ref{th1.2}.
\begin{prop}\label{cor1.3}
Let $\boldsymbol{c}^n = (c^n_1, \ldots , c^n_m)$ be the eigenvector corresponding a simple and $n$-th eigenvalue of the matrix $(h_{ij})$ in (\refeq{f4.1}). If $c^n_j \neq 0$ then $w^n_k$ concentrates at $x^\ast_j$.
\end{prop}
%
%
%
%
\par The following theorem is additional fact considering $c^n_j \neq 0$. 
\begin{thm}\label{th1.4}
Under the same assumption of Theorem \ref{th1.1}, let $w^n_k$, $1 \leq n \leq m$ be the eigenfunctions of (\ref{f1.2}). Then the following (i) and (ii) hold:
\begin{enumerate}[(i)]
  \item $w^1_k$ concentrates at $m$ points $x^\ast_1, \ldots , x^\ast_m \in \mathcal{S}$.
  \item For $m \geq 2$ and each $n= 1, \ldots , m$, $w^n_k$ concentrates at least at two points $x^\ast_i, x^\ast_j \in \mathcal{S}$, $1 \leq i , j  \leq m$, $i \neq j$. 
\end{enumerate}
\end{thm}
Theorems \ref{th1.2} and \ref{th1.4} are also extensions of results in \cite{ggo16}. We will prove those theorems by similar arguments of \cite{ggo16}. In results \cite{oss13} and \cite{ss18} considering the same Liouville-Gel'fand problem (\refeq{f1.1}), we took partial derivatives of the equation (\ref{f1.1}) and used the precise estimates with Taylor expansion of $\log V(x)$ near the critical point of the Hamiltonian $H^m_V$, which cause technical difficulties. However, in this paper, those difficulties do not arise because our argument requires no information of partial derivatives of the equation (\ref{f1.1}) and partial derivatives of $\log V$ except for determining the location of blow-up points. It suggests that there may be similar phenomenon to Theorem \ref{th1.4} under lower regularity of $V$, although we have not seen anything yet.
\par This paper is organized as follows: In section 2, we study the asymptotic behavior of $v_k(x_{j,k})$ more precisely than (\refeq{f6.5}), and we prove Theorem \ref{th1.6} after preparation of related lemmas. In section 3, we study the asymptotic behavior of $w^n_k(x_{j,k}) / \mu^n_k$, and we prove Theorems \ref{th1.2} and \ref{th1.4}.
\section{Analysis for the eigenvalue and proof of Theorem \ref{th1.6}}
In this section, we have a refinement of asymptotic behavior of $n$-th eigenvalue, $1 \leq n \leq m$ for (\refeq{f1.2}) through studying an asymptotic representation for the local maximum of blow-up solution $v_k(x_{j,k})$ more precisely than (\refeq{f6.5}). 
\par First, we prove the following lemma by using Green's formula. Throughout this section, we consider the  assumption of Theorem \ref{th1.1}.
\begin{lem}\label{prop3.1}
It holds that
\begin{align}
&\left\{\frac{1}{\mu^n_k} - v_k(x_{j,k}) \right\} \lambda_k \int_{B_R(x_{j,k})} V e^{v_k} w^n_k \\
&= (8 \pi)^2 \sum_{\substack{1\leq i \leq m \\ i \neq j}} (c^n_i - c^n_j) G(x^\ast_j, x^\ast_i) - 16\pi c^n_j +o(1) \label{f7.1}
\end{align}
for $n = 1, \ldots , m$ and $k \gg 1$. 
\end{lem}
\subsubsection*{Proof of Lemma \ref{prop3.1}}
By (\refeq{f1.1}), (\refeq{f1.2}), (\refeq{f2.1}), (\refeq{f2.2}), (\refeq{f6.7}) and the dominated convergence theorem, we have 
\begin{align}
&\int_{\partial B_R(x_{j,k})} \left\{\frac{\partial v_k}{\partial \nu} \frac{w^n_k}{\mu^n_k} - v_k \frac{\partial}{\partial \nu}\left(\frac{w^n_k}{\mu^n_k} \right) \right\}ds = \int_{B_R(x_{j,k})} \left(\Delta v_k \cdot \frac{w^n_k}{\mu^n_k} - v_k \Delta \frac{w^n_k}{\mu^n_k} \right)ds \\
&= - \frac{1}{\mu^n_k} \int_{B_R(x_{j,k})} \lambda_k V e^{v_k} w^n_k + v_k(x_{j,k}) \int_{B_R(x_{j,k})} \lambda_k V e^{v_k} w^n_k \\
&\quad + \int_{B_R(x_{j,k})} \lambda_k V e^{v_k} \{v_k - v_k(x_{j,k})\} w^n_k \\
&= \left\{ - \frac{1}{\mu^n_k} +  v_k(x_{j,k}) \right\} \int_{B_R(x_{j,k})} \lambda_k V e^{v_k} w^n_k \\
&\quad + \int_{B_{R/\delta_{j,k}}(0)} \tilde{V}_{j,k}(x) e^{\tilde{v}_{j,k}(x)} \tilde{v}_{j,k}(x) \tilde{w}^n_{j,k}(x) dx \\
&= \left\{ - \frac{1}{\mu^n_k} +  v_k(x_{j,k}) \right\} \int_{B_R(x_{j,k})} \lambda_k V e^{v_k} w^n_k + \int_{\mathbb{R}^2} e^U U c^n_j + o(1) \\
&= \left\{ - \frac{1}{\mu^n_k} +  v_k(x_{j,k}) \right\} \int_{B_R(x_{j,k})} \lambda_k V e^{v_k} w^n_k - 16 \pi c^n_j + o(1) \label{f7.2}
\end{align}
where $\nu$ is the outer unit normal vector on $\partial B_R(x_{j,k})$. Concerning the left-hand side of (\refeq{f7.2}), from (\refeq{f0.3}) and (\refeq{f3.4}), it follows that 
\begin{align}
&\int_{\partial B_R(x_{j,k})} \left\{\frac{\partial v_k}{\partial \nu} \frac{w^n_k}{\mu^n_k} - v_k \frac{\partial}{\partial \nu}\left(\frac{w^n_k}{\mu^n_k} \right) \right\}ds \\
&= (8\pi)^2 \sum^{m}_{i=1} \sum^{m}_{h=1} c^n_h \int_{\partial B_R(x^\ast_j)} \left\{\frac{\partial}{\partial \nu} G(x, x^\ast_i) \cdot G(x, x^\ast_h) \right.  \\
& \qquad \qquad \qquad \qquad \qquad \qquad \left. - G(x, x^\ast_i) \frac{\partial}{\partial \nu} G(x, x^\ast_h)  \right\} ds  + o(1) \label{f8.1} \\
&= (8\pi)^2 \sum^{m}_{i=1} \sum^{m}_{h=1} c^n_h I_{i,h} +o(1). 
\end{align}
Here, we have 
\begin{align}
I_{i,h}=
\begin{cases}
0, \quad \text{if} \ i=h, \\
-G(x^\ast_j, x^\ast_h) \delta^j_i + G(x^\ast_j, x^\ast_i) \delta^j_h \label{f8.2}, \quad \text{if} \ i \neq h,
\end{cases}
\end{align}
where $\delta^b_a =1$ if $a=b$, and $\delta^b_a =0$ otherwise.
Relations (\refeq{f8.1}) and (\refeq{f8.2}) indicate 
\begin{align}
&\int_{\partial B_R(x_{j,k})} \left\{\frac{\partial v_k}{\partial \nu} \frac{w^n_k}{\mu^n_k} - v_k \frac{\partial}{\partial \nu}\left(\frac{w^n_k}{\mu^n_k} \right) \right\}ds \\
&= (8\pi)^2 \sum^{m}_{i=1} \sum_{1 \leq h \leq m, \ h \neq i} c^n_h \left\{ -G(x^\ast_j, x^\ast_h) \delta^j_i + G(x^\ast_j, x^\ast_i) \delta^j_h  \right\} +o(1) \\
&= -(8\pi)^2 \sum_{1\leq i \leq m, \ i \neq j}(c^n_i - c^n_j) G(x^\ast_j, x^\ast_i) + o(1). \label{f8.3}
\end{align}
Therefore, we have (\refeq{f7.1}) by (\refeq{f7.2}) and (\refeq{f8.3}). 
\mbox{} \hfill $\Box$
\begin{lem}\label{prop3.2}
Let
\begin{align}
\sigma_{j,k} := \int_{B_R(x_{j,k})} \lambda_k V e^{v_k} = 8 \pi + o(1). \label{f9.2}
\end{align}
Then it holds that
\begin{align}
v_k(x_{j,k}) = &- \frac{\sigma_{j,k}}{\sigma_{j,k} - 4\pi} \left\{\log \lambda_k + \log V(x^\ast_{j}) \right\} + 6 \log 2 \\
&- 8\pi \left\{R(x^\ast_{j}) + \sum_{1\leq i \leq m, \ i \neq j} G(x^\ast_{j}, x^\ast_{i}) \right\} +o(1) \label{f9.1}
\end{align}
where $x_{j,k}$ satisfies (\refeq{f0.12}).
\end{lem}
\subsubsection*{Proof of Lemma \ref{prop3.2}}
By the Green representation formula, it holds that 
\begin{align}
v_k(x_{j,k}) &= \int_{\Omega} G(x_{j,k}, y) \lambda_k V(y) e^{v_k(y)} dy \\
&= - \frac{1}{2\pi} \int_{B_R(x_{j,k})} \log |x_{j,k} - y| \cdot \lambda_k V(y) e^{v_k(y)} dy \\
&\quad +  \int_{B_R(x_{j,k})}  K(x_{j,k}, y) \lambda_k V(y) e^{v_k(y)} dy \\
&\quad + \sum_{1 \leq i \leq m, \ i \neq j} \int_{B_R(x_{j,k})}  G(x_{j,k}, y) \lambda_k V(y) e^{v_k(y)} dy \\
&\quad + \int_{\Omega \setminus \cup^{m}_{j=1} B_R(x_{j,k})}  G(x_{j,k}, y) \lambda_k V(y) e^{v_k(y)} dy  \\
& = I_{1,k}+ I_{2,k}+ I_{3,k}+ I_{4,k}. \label{f9.3}
\end{align}
Concerning $I_{1,k}$ of (\refeq{f9.3}), from (\refeq{f2.1}) we have 
\begin{align}
I_{1,k} &= - \frac{1}{2\pi} \int_{B_{R/\delta_{j,k}}(0)} \log |x_{j,k} - \delta_{j,k} \tilde{y} -x_{j,k}| \\
&\qquad \qquad \qquad \qquad \cdot \lambda_k V(\delta_{j,k} \tilde{y} + x_{j,k}) e^{v_k(\delta_{j,k} \tilde{y} + x_{j,k})} \cdot \delta^2_{j,k} d \tilde{y} \\
&= - \frac{\sigma_{j,k}}{2\pi} \log \delta_{j,k} + \frac{1}{2\pi} \int_{\mathbb{R}^2} \log | \tilde{y}|^{-1} e^{U( \tilde{y})} d \tilde{y} +o(1)\\
&= - \frac{\sigma_{j,k}}{2\pi} \log \delta_{j,k} - 6 \log 2 +o(1) \\
&= \frac{\sigma_{j,k}}{4\pi} \left\{\log \lambda_k +  \log V(x^\ast_{j}) +  v_k(x_{j,k}) \right\} - 6 \log 2 + o(1) \label{f9.4}
\end{align} 
Moreover, 
\begin{align}
I_{2,k} = \int_{B_R(x_{j,k})}  K(x_{j,k}, y) \lambda_k V(y) e^{v_k(y)} dy = 8 \pi R(x^\ast_{j}) + o(1), \label{f11.1}
\end{align}
\begin{align}
I_{3,k} = 8 \pi \sum_{1 \leq i \leq m, \ i \neq j}  G(x^\ast_{j}, x^\ast_{i}) +o(1), \label{f11.2}
\end{align}
and
\begin{align}
I_{4,k} = \int_{\Omega \setminus \cup^{m}_{j=1} B_R(x_{j,k})}  G(x_{j,k}, y) \lambda_k V(y) e^{v_k(y)} dy = o(1) \label{f11.3}
\end{align}
respectively. Therefore, it holds that 
\begin{align}
\left(1- \frac{\sigma_{j,k}}{4\pi} \right) v_k(x_{j,k}) &= \frac{\sigma_{j,k}}{4\pi} \left\{\log \lambda_k + \log V(x^\ast_{j}) \right\} - 6 \log 2 \\
&\quad + 8\pi \left\{R(x^\ast_{j}) + \sum_{1\leq i \leq m, \ i \neq j} G(x^\ast_{j}, x^\ast_{i}) \right\} +o(1) \label{f11.5}
\end{align}
by relations (\refeq{f9.3}), (\refeq{f9.4}), (\refeq{f11.1}), (\refeq{f11.2}), (\refeq{f11.3}). Then we have (\refeq{f9.1}) since $1- \frac{\sigma_{j,k}}{4\pi} \to  -1 \neq 0$. \mbox{} \hfill $\Box$
\begin{lem}\label{prop3.3}
For each $j=1, \ldots , m$ and $0 < \varepsilon < 1/2$, it holds that
\begin{align}
\sigma_{j,k} = 8\pi + o(\lambda_k^{\frac{1}{2}-\varepsilon}). \label{f12.1}
\end{align}
\end{lem}
We note that the sharper estimate 
\begin{align}
\sigma_{j,k} = 8\pi + o(\lambda_k)
\end{align}
is already known, see \cite{cl02}. This is, however, obtained under the condition $V\in C^2(\overline{\Omega})$. We need only $V\in C^1(\overline{\Omega})$ to get Lemma \ref{prop3.3}. 
\subsubsection*{Proof of Lemma \ref{prop3.3}}
We recall the Pohozaev identity in local form, for example, see Proposition 5.5 in \cite{o12},
\begin{align}
2 \int_{B_R(x_{j,k})} (x-x_{j,k}) \cdot \nabla v_k \Delta v_k = R \int_{\partial B_R(x_{j,k})} \left\{ 2 \left(\frac{\partial v_k}{\partial \nu} \right)^2 - |\nabla v_k|^2 \right\}. \label{ef4.3}
\end{align}
Since $v_k$ satisfies (\refeq{f1.1}), it holds that
\begin{align}
&2 \int_{B_R(x_{j,k})} (x-x_{j,k}) \cdot \nabla v_k \Delta v_k = -2 \int_{B_R(x_{j,k})} (x-x_{j,k}) \cdot \nabla v_k \{\lambda_k V(x) e^{v_k}\} \\
&= -2 \int_{B_R(x_{j,k})} (x-x_{j,k}) \cdot \nabla \{\lambda_k V(x) e^{v_k} \} \\
& \quad+ 2 \int_{B_R(x_{j,k})} (x-x_{j,k}) \cdot \nabla V(x) \{\lambda_k e^{v_k} \} = I_{1,k} + I_{2,k}. \label{ef4.5}
\end{align}
Concerning $I_{1,k}$ of (\refeq{ef4.5}), it holds that
\begin{align}
I_{1,k} &= -2 \int_{B_R(x_{j,k})} (x-x_{j,k}) \cdot \nabla \{\lambda_k V(x) e^{v_k} \} \\
&= -2 \int_{\partial B_R(x_{j,k})} \nu \cdot (x-x_{j,k}) \lambda_k V(x) e^{v_k} + 4 \int_{B_R(x_{j,k})} \lambda_k V(x) e^{v_k}\\
&= O(\lambda_k) + 4 \sigma_{j,k} = o(\lambda^{1/2}_k) + 4 \sigma_{j,k}. \label{ef5.1}
\end{align}
Concerning $I_{2,k}$ of (\refeq{ef4.5}), we note 
\begin{align}
\lambda_k e^{v_k} \nabla V(x) = \lambda_k V(x) e^{v_k} \nabla \log V(x) 
\end{align}
and $\|\nabla \log V\|_{\infty} \leq C$ from the assumptions.
Therefore, taking constants $0 < \varepsilon < 1/2$ and $C>0$, it holds that
\begin{align}
|I_{2,k}| &\leq 2 C \int_{B_R(x_{j,k})} \lambda_k V(x) e^{v_k} |x- x_{j,k}| dx \\
&= O(1) \cdot \delta_{j,k} \left\{\int_{\mathbb{R}^2} e^U |x| dx + o(1) \right\} \\
&= O(\delta_{j,k}) = o\left(\lambda^{\frac{1}{2}- \varepsilon}_k \right) \label{ef5.3}
\end{align}
by the estimates (\refeq{f6.4}), (\refeq{f6.7}) and $|x|^\alpha e^U \in L^1(\mathbb{R}^2)$ for $0 < \alpha <2$, where $U=U(x)$ is defined in (\refeq{f2.4}).
\par Concerning the right-hand side of (\refeq{ef4.3}), we recall the formula (29) in \cite{ss18}, that is,
\begin{align}
\nabla v_k(x) = - \frac{\sigma_{j,k}}{2\pi} \cdot \frac{\nu}{R} + \nabla k_{j,k}(x) +o(\lambda^{1/2}_k) \label{ef6.1}
\end{align}
uniformly on $\partial B_R(x_{j,k})$, where
\begin{align}
k_{j,k}(x) = \sigma_{j,k}K(x,x_{j,k}) + \sum_{1 \leq i \leq m, \ i \neq j} \sigma_{j,k} G(x,x_{j,k}) \label{ef6.2}
\end{align}
is harmonic in $B_R(x_{j,k})$.
By (\refeq{ef6.1}) and (\refeq{ef6.2}) it holds that
\begin{align}
&R \int_{\partial B_R(x_{j,k})} \left\{ 2 \left(\frac{\partial v_k}{\partial \nu} \right)^2 - |\nabla v_k|^2 \right\} ds \\
&= R \int_{\partial B_R(x_{j,k})} \left\{ \frac{(\sigma_{j,k})^2}{4\pi^2 R^2} + 2 \left(\frac{\partial k_{j,k}}{\partial \nu} \right)^2 - |\nabla k_{j,k}|^2 - \frac{\sigma_{j,k}}{\pi R} \cdot \frac{\partial k_{j,k}}{\partial \nu} + o(\lambda^{1/2}_k) \right\} ds \\
&= \frac{(\sigma_{j,k})^2}{2 \pi} + R \int_{\partial B_R(x_{j,k})} \left\{ 2 \left(\frac{\partial k_{j,k}}{\partial \nu} \right)^2 - |\nabla k_{j,k}|^2 \right\} ds + o(\lambda^{1/2}_k) \\
&= \frac{(\sigma_{j,k})^2}{2\pi} + o(\lambda^{1/2}_k). \label{ef6.5}
\end{align}
The final step of (\refeq{ef6.5}) comes from (\refeq{ef4.3}) with $k_{j,k}$.
\par From (\refeq{ef4.3}), (\refeq{ef4.5}), (\refeq{ef5.1}), (\refeq{ef5.3}), (\refeq{ef6.5}) it follows that 
\begin{align}
o(\lambda^{1/2}_k) + 4 \sigma_{j,k} + o(\lambda^{\frac{1}{2} - \varepsilon}_k) = \frac{(\sigma_{j,k})^2}{2\pi} + o(\lambda^{1/2}_k) \label{ef7.1}
\end{align}
Relations (\refeq{f9.2}) and (\refeq{ef7.1}) imply (\refeq{f12.1}).
\mbox{} \hfill $\Box$
\par Regarding a representation of constants $d_j$ in (\refeq{f6.4}), we have the next lemma which is also an extension for the case of $V \equiv 1$.
\begin{lem}\label{prop3.4}
For each $j=1, \ldots , m$, it holds that
\begin{align}
d_j = \frac{1}{8} \exp \left\{4\pi R(x^\ast_j) + 4\pi \sum_{1\leq i \leq m, \ i \neq j} G(x^\ast_j, x^\ast_i) + \frac{1}{2} \log V (x^\ast_j)  \right\}. \label{f12.2}
\end{align}
\end{lem}
\subsubsection*{Proof of Lemma \ref{prop3.4}}
From (\refeq{f9.1}), we obtain  
\begin{align}
v_k(x_{j,k}) &= -2 \left\{\log \lambda_k + \log V(x^\ast_j) \right\} + \frac{\sigma_{j,k} - 8\pi}{\sigma_{j,k} -4\pi} \left\{\log \lambda_k + \log V(x^\ast_j) \right\} \\
&\quad - 8 \pi  \left\{R(x^\ast_j) + \sum_{1\leq i \leq m, \ i \neq j} G(x^\ast_j, x^\ast_i) \right\} + 6 \log 2 + o(1) \label{f12.3}
\end{align}
On the other hand, it holds that 
\begin{align}
\frac{\sigma_{j,k} - 8\pi}{\sigma_{j,k} -4\pi} \left\{\log \lambda_k + \log V(x^\ast_j) \right\} = o(1) \label{f12.4}
\end{align}
by (\refeq{f12.1}). Then relations (\refeq{f6.5}), (\refeq{f12.3}) and (\refeq{f12.4}) indicate (\refeq{f12.2}). \mbox{} \hfill $\Box$
\par Here, by using (\refeq{f3.2}) we note 
\begin{align}
\int_{B_R(x_{j,k})} \lambda_k V e^{v_k} w^n_k &= \int_{B_{R/\delta_{j,k}}(0)} \tilde{V}_{j,k}(x) e^{\tilde{v}_{j,k}(x)} \tilde{w}^n_{j,k}(x) dx \\
&= c^n_j \int_{\mathbb{R}^2} e^U + o(1) = 8 \pi c^n_j + o(1). \label{f18.1}
\end{align}
We show the following two lemmas which characterize the vector $\boldsymbol{c}^n$ in (\refeq{f3.1}). 
\begin{lem}\label{prop3.5}
For each $j=1, \ldots , m$ and $\ell=1, \ldots , m$, it holds that
\begin{align}
c^n_\ell \sum^{m}_{i=1} h_{ji} c^n_i = c^n_j \sum^{m}_{i=1} h_{\ell i} c^n_i  \label{f13.2}
\end{align}
where $h_{ij}$ is defined in (\refeq{f4.1}).
\end{lem}
\subsubsection*{Proof of Lemma \ref{prop3.5}}
From (\refeq{f6.5}), (\refeq{f3.4}), (\refeq{f12.2}) and (\refeq{f18.1}), it follows that 
\begin{align}
&\left\{\frac{1}{\mu^n_k} + 2 \log \lambda_k \right\} \int_{B_R(x_{j,k})} \lambda_k V e^{v_k} w^n_k \\
&= \left\{\frac{1}{\mu^n_k} - v_k(x_{j,k}) - 2\log d_j - \log V(x^\ast_j)+ o(1) \right\}  \int_{B_R(x_{j,k})} \lambda_k V e^{v_k} w^n_k \\
&=(8 \pi)^2 \sum_{1\leq i \leq m, \ i \neq j}(c^n_i- c^n_j) G(x^\ast_j, x^\ast_i) - 16 \pi c^n_j \\
&\quad -2 \log V(x^\ast_j) \left\{8\pi c^n_j + o(1) \right\} + 6 \log 2 \cdot \left\{8 \pi c^n_j + o(1) \right\} \\
&\quad - 8\pi \left\{ R(x^\ast_j) +  \sum_{\substack{1\leq i \leq m \\ i \neq j}}G(x^\ast_j, x^\ast_i) \right\} \left\{8\pi c^n_j + o(1) \right\} + o(1) \\
&=-(8\pi)^2 \sum^{m}_{i=1} h_{ji} c^n_i + 16\pi c^n_j (3\log 2 -1) + o(1). \label{f13.1}
\end{align}
By using (\refeq{f13.1}), we obtain 
\begin{align}
0 &= \left\{\frac{1}{\mu^n_k} + 2 \log \lambda_k \right\} \left\{ \int_{B_R(x_{j,k})} \lambda_k V e^{v_k} w^n_k\right\}\left\{ \int_{B_R(x_{\ell, k})} \lambda_k V e^{v_k} w^n_k\right\} \\
&\quad - \left\{\frac{1}{\mu^n_k} + 2 \log \lambda_k \right\} \left\{ \int_{B_R(x_{\ell,k})} \lambda_k V e^{v_k} w^n_k\right\}\left\{ \int_{B_R(x_{j, k})} \lambda_k V e^{v_k} w^n_k\right\} \\
&= \left\{-(8\pi)^2 \sum^{m}_{i=1} h_{ji} c^n_i + 16\pi c^n_j (3\log 2 -1) + o(1) \right\} \left\{8 \pi c^n_{\ell} + o(1) \right\} \\
&\quad - \left\{-(8\pi)^2 \sum^{m}_{i=1} h_{\ell i} c^n_i + 16\pi c^n_{\ell} (3\log 2 -1) + o(1) \right\} \left\{8 \pi c^n_j + o(1) \right\}. \label{f14.1}
\end{align}
Therefore, we have (\refeq{f13.2}) by (\refeq{f14.1}).
\mbox{} \hfill $\Box$
\subsubsection*{Proof of Theorem \ref{th1.6}}
(i) \ This follows from Lemma \ref{prop3.5}, see the proof of Proposition 3.6 in \cite{ggo16}. \\
(ii) \ Choose $j=1, \ldots , m$ satisfying $c^n_j \neq 0$. We obtain 
\begin{align}
\sum^{m}_{i=1} h_{ji} c^n_i = \Lambda^n c^n_j \label{f15.2}
\end{align}
by above (i). Therefore, from (\refeq{f18.1}) and (\refeq{f13.1}) it follows that 
\begin{align}
\frac{1}{\mu^n_k} = -2 \log \lambda_k - 8\pi \Lambda^n + 2(3 \log 2 -1) + o(1) \label{f15.3}
\end{align}
and
\begin{align}
\mu^n_k = \frac{1}{-2 \log \lambda_k + L + o(1)} 
= - \frac{1}{2 \log \lambda_k} - \frac{L}{4} \cdot \frac{1}{(\log \lambda_k)^2} + o\left( \frac{1}{(\log \lambda_k)^2} \right) \label{f15.5}
\end{align}
where $L = - 8 \pi \Lambda^n + 2(3 \log 2 -1)$. Thus we have (\refeq{f5.2}). It holds that $\Lambda^n$ is the $n$-th eigenvalue of $(h_{ij})$ because relations (\refeq{f0.13}) and (\refeq{f5.2}) indicate $\Lambda^1 \leq \cdots \leq \Lambda^m$. Asymptotic behavior (\refeq{f5.2}) is valid without taking a subsequence since $\Lambda^n$ depends on $G(x,y)$, $R(x)$ and $\log V(x)$.
\mbox{} \hfill $\Box$
\section{Analysis for the eigenfunction $w^n_k$ and proofs of Theorems \ref{th1.2} and \ref{th1.4}}
In this section, we obtain Theorems \ref{th1.2} and \ref{th1.4}. 
First, we prepare the following lemma for the study of asymptotic behavior of $w^n_k(x_{j,k}) / \mu^n_k$ along the argument of \cite{ggo16}. 
\begin{lem}\label{prop4.1}
For each $j=1, \ldots , m$ and $n=1, \ldots , m$, it holds that
\begin{align}
\frac{w^n_k(x_{j,k})}{\mu^n_k} 
=& \frac{1}{2\pi} \log \delta^{-1}_{j,k} \cdot \int_{B_R(x_{j,k})} \lambda_k V(y) e^{v_k(y)} w^n_k(y) dy  - 6 c^n_j \log 2 \\
=& \left\{\log \lambda_k + \log V(x^\ast_j) + v_k(x_{j,k}) \right\} \cdot \left\{2 c^n_j +o(1) \right\} - 6 c^n_j \log 2 \\
&+8\pi \left\{c^n_j R(x^\ast_j) + \sum_{1 \leq i \leq m, \ i \neq j} c^n_i G(x^\ast_j , x^\ast_i) \right\} +o(1). \label{f16.1}
\end{align}
\end{lem}
\subsubsection*{Proof of Lemma \ref{prop4.1}}
Similarly to (\refeq{f9.3})-(\refeq{f11.3}), we have (\refeq{f16.1}). \mbox{} \hfill $\Box$
\begin{lem}\label{prop4.3}
For each $j=1, \ldots , m$ and $n=1, \ldots , m$, it holds that
\begin{align}
\int_{B_R(x_{j,k})} \lambda_k  V(x) e^{v_k(x)} \cdot \frac{w^n_k(x) - w^n_k(x_{j,k})}{\mu^n_k} dx = -16 \pi c^n_j +o(1).  \label{f19.1}
\end{align}
\end{lem}
\subsubsection*{Proof of Lemma \ref{prop4.3}}
By (\refeq{f13.1}), it holds that 
\begin{align}
&\frac{1}{\mu^n_k}\int_{B_R(x_{j,k})}\lambda_k V(x) e^{v_k(x)} w^n_k(x) \\
&= -2 \log \lambda_k \cdot \int_{B_R(x_{j,k})} \lambda_k V e^{v_k} w^n_k -(8\pi)^2 \sum^{m}_{i=1} h_{ji} c^n_i + 16\pi c^n_j (3\log 2 -1) \\
& \quad + o(1) \label{f13.1dash}
\end{align}
On the other hand, by (\refeq{f2.1}) and (\refeq{f12.1}), we have 
\begin{align}
&\sigma_{j,k} \log \delta_{j,k}^{-1} = \left\{4\pi + o(\lambda^{\frac{1}{2}- \varepsilon}_k) \right\} \left\{\log \lambda_k + \log V(x_{j,k}) + v_k(x_{j,k}) \right\}. \label{f17.5}
\end{align}
By (\refeq{f9.1}), (\refeq{f18.1}), (\refeq{f16.1}) and (\refeq{f17.5}) it holds that 
\begin{align}
&\frac{1}{\mu^n_k} \int_{B_R(x_{j,k})} \lambda_k V e^{v_k} w^n_k(x_{j,k}) =\sigma_{j,k} \cdot \frac{w^n_k(x_{j,k})}{\mu^n_k} \\
&= -2 \log \lambda_k \cdot \int_{B_R(x_{j,k})} \lambda_k V e^{v_k} w^n_k - (8\pi)^2 \sum^{m}_{i=1}h_{ji} c^n_i  + 48 \pi c^n_j \log 2 \\
& \quad +o(1) \label{f18.2}
\end{align}
where $h_{ij}$ is defined in (\refeq{f4.1}). Therefore, we have (\refeq{f19.1}) by subtracting (\refeq{f18.2}) from (\refeq{f13.1dash}).
\mbox{} \hfill $\Box$
\subsubsection*{Proofs of Theorems \ref{th1.2} and \ref{th1.4}}
See the proofs of Theorems 1.2 and 1.4 of \cite{ggo16} respectively, because it works if the matrix $(h_{ij})$ is substituted by (\refeq{f4.1}). \mbox{} \hfill $\Box$ 
\section*{Acknowledgment}
This work was supported by Japan Society for the Promotion of Science Grant-in-Aid for Scientific Research Numbers JP19H01799, JP20K03675.

%
%
%
\end{document}